\documentclass[12pt, a4paper]{article}

\usepackage[utf8]{inputenc} 
\usepackage[T1]{fontenc} 
\usepackage{lmodern} 
\usepackage{amsmath, amsthm, amssymb, amsfonts} 
\usepackage{geometry} 
\usepackage{fancyhdr} 
\usepackage{titlesec} 
\usepackage{parskip} 
\usepackage{hyperref} 
\usepackage{xcolor} 
\usepackage{graphicx} 
\usepackage{enumitem} 

\fancypagestyle{firstpage}{
    \fancyhf{}
    \fancyhead[L]{\includegraphics[width=1cm]{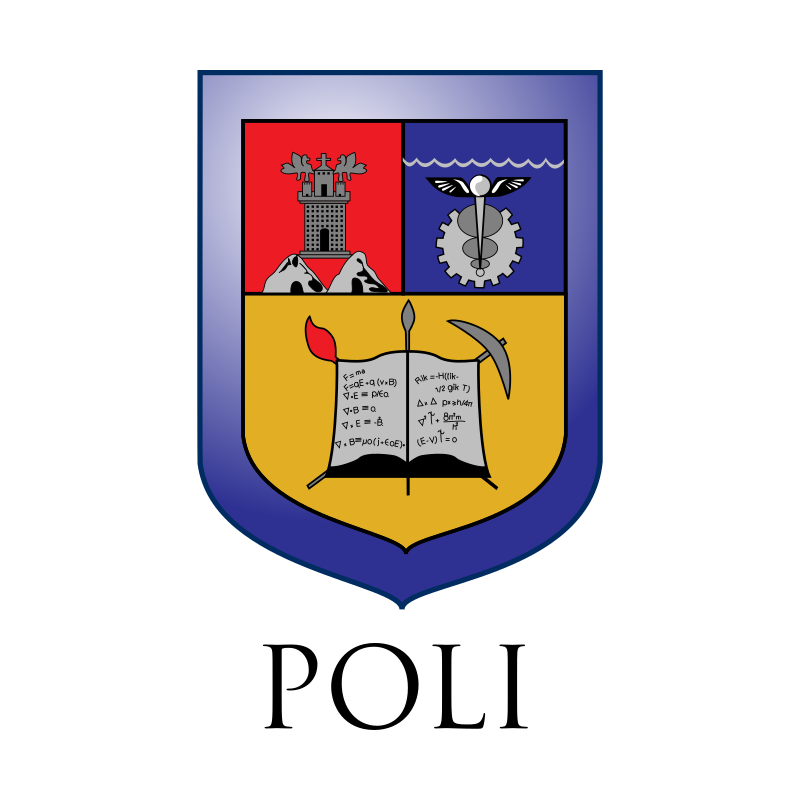}}
    \fancyhead[C]{}
    \fancyhead[R]{\includegraphics[width=1cm]{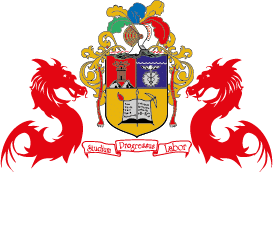}}
    
    \fancyfoot[C]{\thepage}
}

\titleformat{\section}
  {\normalfont\Large\bfseries}{\thesection}{1em}{}
\titleformat{\subsection}
  {\normalfont\large\bfseries}{\thesubsection}{1em}{}
\titleformat{\subsubsection}
  {\normalfont\normalsize\bfseries}{\thesubsubsection}{1em}{}

\newtheoremstyle{mytheoremstyle}
  {\topsep}   
  {\topsep}   
  {\normalfont}  
  {0pt}       
  {\bfseries} 
  {.}         
  { }         
  {\thmname{#1}\thmnumber{ #2}\thmnote{ (#3)}}  

\theoremstyle{mytheoremstyle}
\newtheorem{theorem}{Theorem}
\newtheorem*{vegasthm}{Vega's Theorem} 
\DeclareMathSizes{10}{14}{8}{7} 

\title{Fibonacci sequence and Pythagorean triples in the composition of functions for integer solutions from certain operators}
\author{Pablo José  Vega Esparza\\Banner Code:00336820}
\date{\today}

\begin{document}

\thispagestyle{firstpage}

\maketitle
\subsection*{Abstract}
The following paper synthesizes research where theorems and their respective proofs are postulated based on quadratic equations with special properties given by Phytagorean triples and Fibonacci sequence and an environment in number theory and its applications to calculus. It is important to read the document carefully to build a clear and comprehensive understanding of what is being done.

\section*{}

Integer solutions to an equation are of great importance in number theory. The first question to ask is: What conditions must be met for a quadratic equation to have two real integer solutions?

\begin{align*}
    ax^2 + bx + c = 0
\end{align*}

The following conditions are proposed to find two integer and real roots:

\textbf{formula(1)}:

\[
\boxed{
\begin{aligned}
    &a : \text{Leg of Pythagorean triple} \\\\
    &b = 2a\psi , \quad \psi : \text{Hypotenuse of Pythagorean triple} \\\\
    &c = a^3 \\
\end{aligned}
}
\]

\begin{align*}
    ax^2 + bx + c = 0 ,  \quad
    x\in \mathbb{Z}
\end{align*}

\newpage
\begin{minipage}{\textwidth}
\begin{flushleft}
\textbf{Proof:}
\end{flushleft}

Let \( a, b, c \) belong to integers.

If \( x = \frac{{-b \pm \sqrt{b^2 - 4ac}}}{{2a}} \),

Then, \( x = -\frac{b}{2a} \pm \frac{\sqrt{b^2 - 4ac}}{2a} \).

With \( b \), it is deduced that it must have among its factors \( a \), along with \( 2 \), so \( b = 2 \cdot a \cdot \psi \), where \( \psi \) is the factor by which \( b \) will be a multiple of \( a \) along with 2, and therefore, proportional.

\[
\implies x = -\frac{2 \cdot \psi \cdot a}{2 \cdot a} \pm \frac{1}{2} \sqrt{\frac{(2 \cdot \psi \cdot a)^2 - 2^2 \cdot a \cdot c}{a^2}}
\]

\[
\implies x = -\frac{2 \cdot \psi \cdot a}{2 \cdot a} \pm \sqrt{\psi^2 -\frac{ c}{a}}\]

\[
\implies x = -\psi \pm \sqrt{\psi^2 -\frac{ c}{a}}
\]

 If \( c = a^3 \), notice that it forms a perfect square of \( a \), therefore we get

\begin{equation}
    x = -\psi \pm \sqrt{\psi^2 - a^2} \quad \textbf{formula (2)}
    \label{eq:formula2}
\end{equation}
 For the expression to make sense \( \sqrt{\psi^2 - a^2} > 0 \) with \( \psi^2 > a^2 \). Note that if we choose the factor \( \psi \) of \( b \) as the hypotenuse of a Pythagorean triple, along with \( a \) as one of the legs belonging to a Pythagorean triple, we have obtained two integer summands from which two integer solutions are obtained.

\end{minipage}

\[
x_1 \in \mathbb{Z}, \quad x_2 \in \mathbb{Z} \hfill \quad\blacksquare 
\]
\begin{minipage}{\textwidth}
\framebox[\textwidth]{%
    \begin{minipage}{\dimexpr\textwidth-2\fboxsep-2\fboxrule}
        \[
        \begin{aligned}
            & a : \text{Leg of Pythagorean triple} \\
            & b = 2a \cdot \psi , \quad \psi : \text{Hypotenuse of Pythagorean triple} \\
            & c = a^3
        \end{aligned}
        \]
    \end{minipage}%
}\end{minipage}

\section*{Fibonacci Sequence in Pythagorean Triples}

There are families of Pythagorean triples composed of elements from the Fibonacci sequence, which yield special results given the nature of the sequence.

\begin{theorem}:
Let $\varphi_i, \varphi_{i+1}, \varphi_{i+2}, \varphi_{i+3} \in \varphi_n$ be terms of the Fibonacci sequence. Then:

\textbf{formula(3)}
\[
\begin{cases}
\alpha = \varphi_i \cdot \varphi_{i+3} \\
\beta = 2 \cdot \varphi_{i+1} \cdot \varphi_{i+2} \\
\gamma = \varphi_{i+1}^2 + \varphi_{i+2}^2
\end{cases}
\]
satisfy $\alpha^2 + \beta^2 = \gamma^2$, where $\alpha, \beta, \gamma \in \mathbb{Z}$.
\end{theorem}

\textbf{Proof:}
The Fibonacci sequence is defined recursively as follows:
\[
\begin{cases}
\varphi_0 = 0 \\
\varphi_1 = 1 \\
\varphi_n = \varphi_{n-1} + \varphi_{n-2}, \quad \text{for } n \geq 2
\end{cases}
\]

where \(\varphi_n\) represents the \(n\)-th term of the Fibonacci sequence.

The recursive definition implies that $ \varphi_{i+2} = \varphi_{i+1} + \varphi_i \quad \text{and} \quad \varphi_{i+3} = \varphi_{i+2} + \varphi_{i+1}$.

Thus, if we define \fbox{$\alpha = \varphi_i \cdot \varphi_{i+3}$}, we get $\alpha=\varphi_i (\varphi_{i+2} + \varphi_{i+1})$. For \fbox{$\beta = 2 \cdot \varphi_{i+1} \cdot \varphi_{i+2}$}, we have $\beta=2\cdot\varphi_{i+1}\cdot(\varphi_{i}+\varphi_{i+1})$. Finally, for \fbox{$\gamma = \varphi_{i+1}^2 + \varphi_{i+2}^2$}, we obtain $\gamma=\varphi_{i+1}^2 + (\varphi_{i}+\varphi_{i+1})^2$.

Assume that $\alpha$, $\beta$, $\gamma$ do not satisfy $mcd(\varphi_i, \varphi_{i+1}, \varphi_{i+2}, \varphi_{i+3}) \neq 1$.

\[
\alpha^2 + \beta^2 \neq \gamma^2 \implies (\varphi_i (\varphi_{i+2} + \varphi_{i+1}))^2 + (2\cdot\varphi_{i+1}\cdot(\varphi_{i}+\varphi_{i+1}))^2 \neq (\varphi_{i+1}^2 + (\varphi_{i}+\varphi_{i+1})^2)^2
\]

Expanding the expressions, we get:

\[
\left( \scriptstyle{\varphi_i^4 + 4\varphi_i^3\varphi_{i+1} + 4\varphi_i^2\varphi_{i+1}^2} \right) + \left( \scriptstyle{4\varphi_i^2\varphi_{i+1}^2 + 8\varphi_i\varphi_{i+1}^3 + 4\varphi_{i+1}^4} \right) \\
\neq \left( \scriptstyle{\varphi_i^4 + 4\varphi_i^2\varphi_{i+1}^2 + 4\varphi_{i+1}^4 + 4\varphi_i^3\varphi_{i+1} + 4\varphi_i^2\varphi_{i+1}^2 + 8\varphi_i\varphi_{i+1}^3} \right)
\]

Simplifying further:

\[
\varphi_i^4 + 4\varphi_{i+1}^4 + 4\varphi_i^3\varphi_{i+1} + 8\varphi_i\varphi_{i+1}^3 + 8\varphi_i^2\varphi_{i+1}^2 \\
\neq \varphi_i^4 + 4\varphi_{i+1}^4 + 4\varphi_i^3\varphi_{i+1} + 8\varphi_i\varphi_{i+1}^3 + 8\varphi_i^2\varphi_{i+1}^2
\]

This leads to a contradiction, so the assumption must be false. Therefore, 
\[
\alpha^2 + \beta^2 = \gamma^2
\]
with $mcd(\varphi_i, \varphi_{i+1}, \varphi_{i+2}, \varphi_{i+3}) = 1$.

Thus, 
\[
(\varphi_i \cdot \varphi_{i+3})^2 + (2 \cdot \varphi_{i+1} \cdot \varphi_{i+2})^2 = (\varphi_{i+1}^2 + \varphi_{i+2}^2)^2
\]
forms a Pythagorean triple $\blacksquare$.

The connection between the formula determined in the first section and the formula for generating Pythagorean triples from the Fibonacci sequence is significant, as we can write formula (1) in the following way:

\begin{enumerate}
    \item $a = \alpha$, $b = 2 \cdot \alpha \cdot \gamma$, $c = \alpha^3$
    \item $a = \beta$, $b = 2 \cdot \beta \cdot \gamma$, $c = \beta^3$
\end{enumerate}

Note that there are two expressions since the formula works with one of the two legs of a Pythagorean triple, with \(\alpha, \beta\) being legs of the same triple.

Similarly, we can rewrite this expression according to formula (1), given the previously proven equivalence with terms of the Fibonacci sequence:

\begin{enumerate}
    \item \(a = \varphi_i \cdot \varphi_{i+3}\), \(b = 2 \cdot (\varphi_i \cdot \varphi_{i+3}) \cdot (\varphi_{i+1}^2 + \varphi_{i+2}^2)\), \(c = (\varphi_i \cdot \varphi_{i+3})^3\)
    \item \(a = 2 \cdot (\varphi_{i+1} \cdot \varphi_{i+2})\), \(b = 2^2 \cdot (\varphi_{i+1} \cdot \varphi_{i+2}) \cdot (\varphi_{i+1}^2 + \varphi_{i+2}^2)\), \(c = 2^3 \cdot (\varphi_{i+1} \cdot \varphi_{i+2})^3\)
\end{enumerate}

In the next section, sequences characterizing the given formula will be presented due to their connection with the Pythagorean triples generated by the Fibonacci sequence.

\section*{Primitive Pythagorean Triples and Consecutive Patterns}

\textbf{Definition 1:} A Pythagorean triple \((a, b, c)\) is called primitive if its components are relatively prime, that is, if \(mcd(a, b, c) = 1\).

Primitive Pythagorean triples serve as bases for others, considering triples that differ by a positive integer from their base.

\begin{theorem}$:$ Let \(a, b, c \in \mathbb{Z}\) such that \(a^2 + b^2 = c^2 \quad \implies (ka)^2 + (kb)^2 = (kc)^2\), where $k \in \mathbb{Z}$ forms a non-primitive Pythagorean triple.
\end{theorem}
\textbf{Demonstration}:

The demonstration is trivial since it results merely from scaling the expression with an integer factor. It is directly related to the similarity of triangles, which maintain the same properties. \(\blacksquare\)

Once the primitive and non-primitive Pythagorean triples are defined, we proceed to "experiment with what we have."

Let \( f:\mathbb{R}\rightarrow\mathbb{R} \) be defined as:
\[ f(x)= (3+3n)x^2 + 2(3+3n)(5+5n)x + (3+3n)^3 \]
where \( n \in \mathbb{N} \).

Let \( g:\mathbb{R}\rightarrow\mathbb{R} \) be defined as:
\[ g(x)= (4+4n)x^2 + 2(4+4n)(5+5n)x + (4+4n)^3 \]
where \( n \in \mathbb{N} \).

Taking 4 terms of the Fibonacci sequence:
\[
\begin{cases}
\varphi_1 = 1\\
\varphi_2 = 1\\
\varphi_3 = 2\\
\varphi_4 = 3\\
\end{cases}
\]

Determining the values that will form the primitive triple from the 4 consecutive terms of the sequence using formula (3):

\begin{align*}
&\alpha  = \varphi_1 \cdot \varphi_4 = 1 \cdot 3 = 3 \qquad \\\\
&\beta  = 2 \cdot \varphi_3 \cdot \varphi_2 = 2 \cdot 2 \cdot 1 = 4 \qquad \\\\
&\gamma  = \varphi_2^2 + \varphi_3^2 = 1^2 + 2^2 = 5 \qquad
\end{align*}

\textbf{From this primitive triple, we can generate a family of Pythagorean triples by scaling the base triple (3, 4, 5) with an integer. For each integer \( n \) chosen, different sets \((a,b,c)\) of these triples are taken, but they belong to the same family.}

Note that the following expression is true according to Theorem 1.2:
\[
(3+3n)^2 + (4+4n)^2 = (5+5n)^2
\]

If the reader is still not convinced, here is the proof.

\textbf{Proof:}

Consider the following equality:
\[
(3+3n)^2 + (4+4n)^2 = (5+5n)^2
\]

Factoring out \((n+1)^2\):
\[
(n+1)^2 (3^2 + 4^2) = 5^2 (n+1)^2
\implies 3^2 + 4^2 = 5^2
\]

Thus, it is demonstrated. \(\blacksquare\)

Now, note that the functions \( f \) and \( g \) comply with the values to form the coefficients proposed in formula (1):

\[
\boxed{
\begin{aligned}
    &a : \text{Leg of Pythagorean triple} \\
    &b = 2a\psi , \quad \psi : \text{Hypotenuse of Pythagorean triple} \\
    &c = a^3 \\
\end{aligned}
}
\]

We continue with the expression for \( x \) to obtain integer roots:
\[
x = -\psi \pm \sqrt{\psi^2 - a^2}
\]

Replacing the corresponding values in the expression for \( f \) and \( g \):

For \( f \):
\[ x = -(5+5n) \pm \sqrt{(5+5n)^2 - (3+3n)^2} \]

For \( g \):
\[ x = -(4+4n) \pm \sqrt{(5+5n)^2 - (4+4n)^2} \]

We can rewrite the roots for both \( f \) and \( g \):

For \( f \):
\[ x = -(5+5n) \pm (4+4n) \]

For \( g \):
\[ x = -(5+5n) \pm (3+3n) \]

This implies that the roots of \( f \) are given by \( x_1 = -(n+1) \) and \( x_2 = -9(n+1) \).

The roots of \( g \) are given by \( x_1 = -2(n+1) \) and \( x_2 = -8(n+1) \).

The two roots for this family of Pythagorean triples result from subtracting and adding the hypotenuse with the leg that was not chosen for the coefficients of the function, and they all follow a pattern.

Continuing with the roots of the derivative with respect to each function, which will give us the point in \( \text{Dom}(f) \), \( \text{Dom}(g) \) for which the maximum or minimum of \( \text{Ran}(f) \), \( \text{Ran}(g) \) is found:

With respect to \( f \):
\[
f'(x) = 2(3 + 3n)x + 2(3 + 3n)(5 + 5n)
\]
Setting the derivative to zero and solving for \( x \):

\[
2(3 + 3n)x + 2(3 + 3n)(5 + 5n) = 0
\]

Dividing both sides by \( 2(3 + 3n) \):

\[
x + (5 + 5n) = 0
\]

Subtracting \( (5 + 5n) \) from both sides:

\[
x = - (5 + 5n)
\]

Thus, the root of the derivative function is \( x = -5(n + 1) \).

With respect to \( g \):
\[
g'(x) = 2(4+4n)x + 2(4+4n)(5+5n)
\]

Setting the derivative to zero and solving for \( x \):
\[
2(4+4n)x + 2(4+4n)(5+5n) = 0\]

Dividing both sides by \( 2(4+4n) \):
\[
x + (5+5n) = 0
\]

Subtracting \( (5+5n) \) from both sides:
\[
x = - (5 + 5n)
\]

Since the roots for each function are integers, their maximum or minimum will also be integers.

The critical points for each function, regardless of how the coefficients are composed following formula (1), will be the same. Why do we talk about maxima when clearly the functions \( f \) and \( g \) have minima? If we consider formula (1) such that:
\[
\boxed{
\begin{aligned}
    &\textbf{formula 1.1} \\
    &a : -\text{Leg of Pythagorean triple} \\
    &b = -2a\psi , \quad \psi : \text{Hypotenuse of Pythagorean triple} \\
    &c = a^3 \\
\end{aligned}
}
\]

There will get exactly the same results for \( f \) and \( g \) but with opposite signs, both for the roots and the critical points.

The reader may think that the connection between the triples and the given values for the quadratic equation coefficients is not closely related; however, this is incorrect, as the value corresponding in the range of both \( f \) and \( g \) is influenced by the nature of the Pythagorean triples.

With respect to \( f \):

\[
(3 + 3n) \left(-(5 + 5n)\right)^2 + 2(3 + 3n)(5 + 5n) \left(-(5 + 5n)\right) + (3 + 3n)^3 = -4 \cdot 12(n + 1)^3
\]

With respect to \( g \):

\[
(4 + 4n) \left(-(5 + 5n)\right)^2 + 2(4 + 4n)(5 + 5n) \left(-(5 + 5n)\right) + (4 + 4n)^3 = -3 \cdot 12(n+1)^3
\]

It can be noted that the values corresponding to the maximum in \(\text{Ran}(f)\) and \(\text{Ran}(g)\) contain within their factors the leg of the primitive triple \((3, 4, 5)\) that was not chosen for the composition of the function's coefficients, along with a factor which is 12, resulting from the sum of the three sides that make up the right triangle forming the triple. (The calculation is omitted to simplify the document, as it is extensive).

What about the integral from root to root? Could the reader imagine that the integral from one root to another follows a sequence and is also integer for every set \((a,b,c)\) of triples in this family? It turns out to be true.

With respect to \( f: \)

\begin{equation}
   \left| \int_{-9(n+1)}^{-(n+1)} \left( (3 + 3n)x^2 + \left(2(3 + 3n)(5+5n)x + (3+ 3n)^{3} \right) \right) \, dx \right| = 4^4(n+1)^4
\end{equation}

The proof is left as an exercise for the reader (the omission of the calculation is made to simplify the document, as the calculation is exhausting).

Now with respect to \( g: \)

\begin{equation}
     \left| \int_{-8(n+1)}^{-2(n+1)} \left((4 + 4n)x^2 + (2(4 + 4n)(5 + 5n))x +(4+ 4n)^{3}\right) , dx \right| =12^2(n+1)^4
\end{equation}

As mentioned earlier, if the factors for formula (1) are taken as negative due to symmetry, the same roots, critical point, and integral with opposite sign are maintained; in this case, positive. \( \blacksquare \).

\begin{figure}[h]
    \centering
    \includegraphics[scale=0.44]{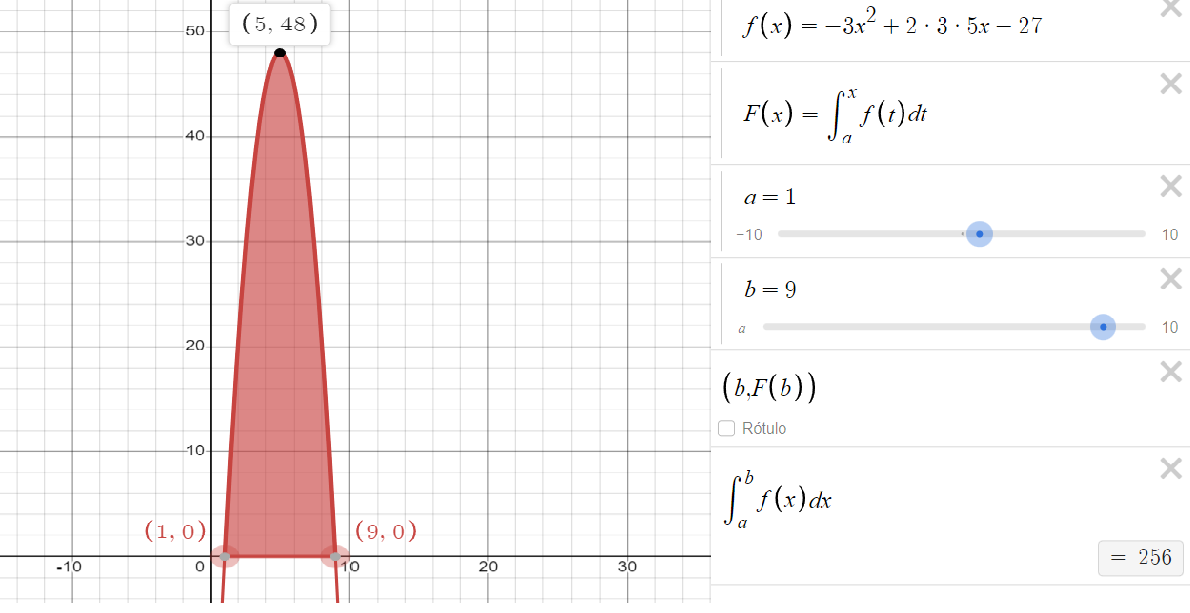}
    \caption{Representation of the integral with respect to the function \( f \), taking \( n = 0 \), considering the symmetry given by formula 1.1.}
    \label{fig:integral}
\end{figure}
\begin{figure}[h]
    \centering
    \includegraphics[scale=0.44]{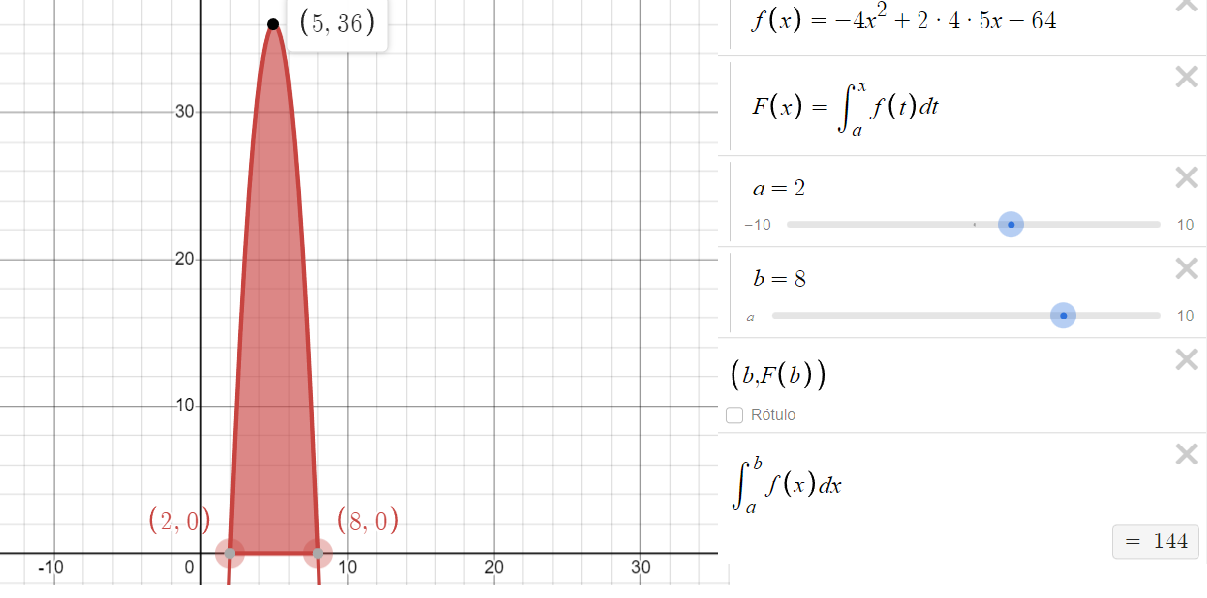}
    \caption{Representation of the integral with respect to the function \( g \), taking \( n = 0 \), considering the symmetry given by formula 1.1.}
    \label{fig:integral}
\end{figure}

The results obtained above sound amazing, right? The reader could affirm that mathematics
couldn't be more generous than he already has been. Family of functions that have the same properties.
It would be easy to search through the factors of a quadratic equation and find the fact that it satisfies
formula (1) and know its roots, minimum or maximum, and value of the integral from root to root, which coincidentally is a perfect square for each Pythagorean triple selected in this family. The patterns found for different results depend on how the functions are composed with the terms of the Fibonacci sequence.\\
Could it be that these properties hold for every Pythagorean triple generated by the Fibonacci sequence?
together with formula (1)?

\section*{Generalization and Theorem to be Presented}

\textbf{Note:} \textbf{It is important for the reader to have become sufficiently familiar with formulas (1), (2), and (3). If not, please refer back to them on the previous pages.}

As already defined based on formula (1) and formula (3),

$*(1).$
$a= \varphi_i \cdot \varphi_{i+3}$, $b=2\cdot(\varphi_i \cdot \varphi_{i+3}) \cdot (\varphi_{i+1}^2 + \varphi_{i+2}^2)$, $c= (\varphi_i \cdot \varphi_{i+3})^3$
\\
$\implies y_1=(\varphi_i \cdot \varphi_{i+3})x^2+2\cdot(\varphi_i \cdot \varphi_{i+3}) \cdot (\varphi_{i+1}^2 + \varphi_{i+2}^2)x+
(\varphi_i \cdot \varphi_{i+3})^3
$\\

$*(2).$
$a= 2 \cdot (\varphi_{i+1} \cdot \varphi_{i+2}) $, $b=2^2\cdot(  (\varphi_{i+1} \cdot \varphi_{i+2})) \cdot (\varphi_{i+1}^2 + \varphi_{i+2}^2)$, $c= 2^3 \cdot( (\varphi_{i+1} \cdot \varphi_{i+2}))^3$\\
$\implies y_2=2 \cdot (\varphi_{i+1} \cdot \varphi_{i+2})x^2+2^2\cdot(  (\varphi_{i+1} \cdot \varphi_{i+2})) \cdot (\varphi_{i+1}^2 + \varphi_{i+2}^2)x+2^3 \cdot (\varphi_{i+1} \cdot \varphi_{i+2})^3$
$\\$

Being \(a\), \(b\), \(c\) coefficients of

\begin{align*}
     y=ax^2 + bx + c ,  \quad
\end{align*}

with formula (2) one can obtain its roots. 

With respect to *(1)
\begin{align*}x=
    -(\varphi_{i+1}^2 + \varphi_{2}^2) \pm \sqrt{(\varphi_{i+1}^2 + \varphi_{i+2}^2)^2 - (\varphi_{i}\cdot\varphi_{i+3})^2} \\\implies
    x=-(\varphi_{i+1}^2 + \varphi_{i+2}^2) \pm(2\cdot \varphi_{i+1} \cdot \varphi_{i+2}) \\\\
    x_1=-(\varphi_{i+1} - \varphi_{i+2})^2\\
    x_2=-(\varphi_{i+1} + \varphi_{i+2})^2
\end{align*}
Whit respect on *(2)

\begin{align*}
      x=-(\varphi_{i+1}^2 + \varphi_{i+2}^2) \pm \sqrt{(\varphi_{i+1}^2 + \varphi_{i+2}^2)^2 - (2\cdot\varphi_{i}\cdot\varphi_{i+2})^2}\\
      \implies x=(\varphi_{i+1}^2 + \varphi_{i+2}^2)\pm(\varphi_{i}\cdot\varphi_{i+3})\\\\
      x_1=-(\varphi_{i+1}^2 + \varphi_{i+2}^2)+(\varphi_{i}\cdot\varphi_{i+3})\\
      x_2=-(\varphi_{i+1}^2 + \varphi_{i+2}^2)-(\varphi_{i}\cdot\varphi_{i+3})
\end{align*}

Once the logic behind the limits is understood, the theorem is presented.

\begin{vegasthm}{3}\\
Let $f:\mathbb{R}\rightarrow\mathbb{R}:f(x)=(\varphi_i \cdot \varphi_{i+3})x^2+2\cdot(\varphi_i \cdot \varphi_{i+3}) \cdot (\varphi_{i+1}^2 + \varphi_{i+2}^2)x+
(\varphi_i \cdot \varphi_{i+3})^3 \\\\
$Let $g:\mathbb{R}\rightarrow\mathbb{R}:g(x)=2 \cdot (\varphi_{i+1} \cdot \varphi_{i+2})x^2+2^2\cdot(  (\varphi_{i+1} \cdot \varphi_{i+2})) \cdot (\varphi_{i+1}^2 + \varphi_{i+2}^2)x+2^3 \cdot (\varphi_{i+1} \cdot \varphi_{i+2})^3$
\\

With $\varphi_i, \varphi_{i+1}, \varphi_{i+2}, \varphi_{i+3}
\in \varphi_n $ (Fibonacci sequence)
\\

$\theta_1=-(\varphi_{i+1} - \varphi_{i+2})^2$, $\theta_2=-(\varphi_{i+1} +\varphi_{i+2})^2$, roots of $f$
\\
\\

$\phi_1=-(\varphi_{i+1}^2 + \varphi_{i+2}^2)+(\varphi_{i}\cdot\varphi_{i+3})$, \quad$\phi_2=-(\varphi_{i+1}^2 + \varphi_{i+2}^2)-(\varphi_{i}\cdot\varphi_{i+3})$, roots of $g$
\\
\\
\begin{align*}
    \implies \int\limits_{\theta_2}^{\theta_1} f(x) \, dx  \quad, \quad \int\limits_{\phi_2}^{\phi_1} g(x) \, dx\quad  \in  \mathbb{Z}
\end{align*}
\end{vegasthm}

\textbf{Proof:}
\\

(1). Let $\varphi_n$ be the Fibonacci sequence, defined recursively as $F_0 = 0$, $F_1 = 1$, and $F_n = F_{n-1} + F_{n-2}$ for $n \geq 2$. The general formula for the $n$th term of the Fibonacci sequence, denoted as $\varphi_n$, can be expressed as:

\begin{equation*}
\varphi_n = \frac{1}{\sqrt{5}} \left( \left(\frac{1+\sqrt{5}}{2}\right)^n - \left(\frac{1-\sqrt{5}}{2}\right)^n \right)
\end{equation*}
\\

The hypothesis to be stated is denoted as
\\
\begin{align*}
    \varphi_{4n} = \frac{1}{\sqrt{5}} \left( \left(\frac{1+\sqrt{5}}{2}\right)^{4n} - \left(\frac{1-\sqrt{5}}{2}\right)^{4n} \right)\equiv0 \quad mod(3).
    \forall \quad n \in \mathbb{N}
\end{align*}
\\
\textbf{Proof(1)}

By the principle of induction, let's denote the base case

\begin{align*}
    &p(1)=\varphi_{4n} = \frac{1}{\sqrt{5}} \left( \left(\frac{1+\sqrt{5}}{2}\right)^{4n} - \left(\frac{1-\sqrt{5}}{2}\right)^{4n} \right) \\
    &\implies \varphi_{4} = \frac{1}{\sqrt{5}} \left( \left(\frac{1+\sqrt{5}}{2}\right)^{4} - \left(\frac{1-\sqrt{5}}{2}\right)^{4} \right) \\
    &\implies = \frac{1}{\sqrt{5}\cdot 2^{4}} \left( \left({1+\sqrt{5}}\right)^{4} - \left({1-\sqrt{5}}\right)^{4} \right) \\
    &= \frac{1}{\sqrt{5}\cdot 2^{4}}\cdot \left( \left(1+\sqrt{5}\right)^{2} + \left(1-\sqrt{5}\right)^{2} \right) \cdot \left( \left(1+\sqrt{5}\right)^{2} - \left(1-\sqrt{5}\right)^{2} \right)\\
    &= \frac{1}{\sqrt{5}\cdot 2^{4}}\cdot  \left(\left(1+2 \sqrt{5}+5\right) + \left(1-2\sqrt{5}+5 \right)\right) \cdot \left(\left(1+2 \sqrt{5}+5\right) - \left(1-2\sqrt{5}+5 \right)\right)
\end{align*}
\begin{align*}
     \\
    \implies \frac{1}{\sqrt{5}}\cdot\frac{4 \sqrt{5}\cdot 12}{ 4^{2}}= 3
\end{align*}

Clearly

\begin{align*}
    3\equiv 0 \quad mod(3)
\end{align*}
hence the base case is proven. Now let's assume

\begin{align*}
    \varphi_{4n} = \frac{1}{\sqrt{5}} \left( \left(\frac{1+\sqrt{5}}{2}\right)^{4n} - \left(\frac{1-\sqrt{5}}{2}\right)^{4n} \right)\equiv0 \quad mod(3).
    \forall \quad n \in \mathbb{N}
\end{align*}
is true
$\implies$
\begin{align*}
    \varphi_{4(n+1)} = \frac{1}{\sqrt{5}} \left( \left(\frac{1+\sqrt{5}}{2}\right)^{4(n+1)} - \left(\frac{1-\sqrt{5}}{2}\right)^{4(n+1)} \right)\equiv0 \quad mod(3).
    \forall \quad n \in \mathbb{N}
\end{align*}
should be true. 
\\
Algebraically developing the expression for $\varphi_{4(n+1)}$

\begin{align*}
     \frac{1}{\sqrt{5}} \left( \left(\frac{1+\sqrt{5}}{2}\right)^{4n} \cdot \left(\frac{1+\sqrt{5}}{2}\right)^{4} - \left(\frac{1-\sqrt{5}}{2}\right)^{4n} \cdot \left(\frac{1-\sqrt{5}}{2}\right)^{4} \right) \\
\end{align*}

\begin{align*}
    \frac{1}{\sqrt{5}} \left( \frac{1+\sqrt{5}}{2} \right)^{4n} \cdot \left( \frac{1+\sqrt{5}}{2} \right)^{4} 
    \quad - \frac{1}{\sqrt{5}} \left( \frac{1-\sqrt{5}}{2} \right)^{4n} \cdot \left( \frac{1-\sqrt{5}}{2} \right)^{4} 
\end{align*}

Without modifying the expression, we proceed with

\begin{align*}
&\frac{1}{\sqrt{5}}\left( \frac{1+\sqrt{5}}{2} \right)^{4n} \cdot \left( \frac{1+\sqrt{5}}{2} \right)^{4} - \frac{1}{\sqrt{5}}\left( \frac{1-\sqrt{5}}{2} \right)^{4n} \cdot \left( \frac{1+\sqrt{5}}{2} \right)^{4} \\
&+ \frac{1}{\sqrt{5}}\left( \frac{1-\sqrt{5}}{2} \right)^{4n} \cdot \left( \frac{1+\sqrt{5}}{2} \right)^{4} - \frac{1}{\sqrt{5}}\left( \frac{1-\sqrt{5}}{2} \right)^{4n} \cdot  \left( \frac{1-\sqrt{5}}{2} \right)^{4} 
\\
\end{align*}
\begin{align*}
= \frac{1}{\sqrt{5}}\cdot \left(\frac{1+\sqrt{5}}{2}\right)^{4}\left(\left(\frac{1+\sqrt{5}}{2} \right)^{4n}  -\left( \frac{1-\sqrt{5}}{2} \right)^{4n}\right)
\quad +\frac{1}{\sqrt{5}} \left( \frac{1 - \sqrt{5}}{2} \right)^{4n} \cdot
\end{align*}
\begin{align*}
&\left( \left(1 +\frac{\sqrt{5}}{2} \right)^4 -\left( \frac{1 - \sqrt{5}}{2} \right)^4 \right)
\end{align*}

By the already proven base case, the expression is rewritten.

\begin{align*}
= \frac{1}{\sqrt{5}}\cdot \left(\frac{1+\sqrt{5}}{2}\right)^{4}\left(\left(\frac{1+\sqrt{5}}{2} \right)^{4n}  -\left( \frac{1-\sqrt{5}}{2} \right)^{4n}\right)
\quad +3 \left( \frac{1 - \sqrt{5}}{2} \right)^{4n} 
\end{align*}

$\implies$

By inductive hypothesis, the assertion holds for $\varphi_{4(n+1)}$
\\
\begin{align*}
\implies
    \varphi_{4n} = \frac{1}{\sqrt{5}} \left( \left(\frac{1+\sqrt{5}}{2}\right)^{4n} - \left(\frac{1-\sqrt{5}}{2}\right)^{4n} \right)\equiv0 \quad mod(3).
    \forall \quad n \in \mathbb{N} \blacksquare
\end{align*}

From what we have proven above, we can deduce that 
$\exists\quad X: X=\{\varphi_i, \varphi_{i+1},\varphi_{i+2},\varphi_{i+3}\} $ with $X\subset\varphi_n$, $ \quad\exists\quad\varphi_k \in X: \varphi_k \equiv0\quad mod(3) \quad *(1)$
\\
Continuing with the proof, regarding $f$
\begin{align*}
\int\limits_{\theta_2}^{\theta_1} f(x) \, dx  \quad=\int\limits_{\theta_2}^{\theta_1}(\varphi_i \cdot \varphi_{i+3})x^2+2\cdot(\varphi_i \cdot \varphi_{i+3}) \cdot (\varphi_{i+1}^2 + \varphi_{i+2}^2)x+
(\varphi_i \cdot \varphi_{i+3})^3 \, dx
\end{align*}

Let's denote for ease of notation
\begin{align*}
    \int\limits_{\theta_2}^{\theta_1}(\varphi_i \cdot \varphi_{i+3})x^2 \,dx+  \int\limits_{\theta_2}^{\theta_1}2\cdot(\varphi_i \cdot \varphi_{i+3}) \cdot (\varphi_{i+1}^2 + \varphi_{i+2}^2)x \, dx+  \int\limits_{\theta_2}^{\theta_1}
(\varphi_i \cdot \varphi_{i+3})^3 \,dx=P1+P2+P3
\end{align*}
\begin{align*}
P1+ \left[  (\varphi_i \cdot \varphi_{i+3}) \cdot (\varphi_{i+1}^2 + \varphi_{i+2}^2)x^2 \right]_{\theta_2}^{\theta_1}
 +\left[ (\varphi_i \cdot \varphi_{i+3}) x \right]_{\theta_2}^{\theta_1}
\end{align*}
Clearly, both $P2$ and $P3$ are integers, integrating $P1$ \\
$\implies$
\begin{align*}
\int\limits_{\theta_2}^{\theta_1} (\varphi_i \cdot \varphi_{i+3})x^2 \, dx  =  [ \frac{1}{3} (\varphi_i \cdot \varphi_{i+3})x^3 ]_{\theta_2}^{\theta_1} 
\end{align*}

\begin{align*}
    \frac{1}{3} (\varphi_i \cdot \varphi_{i+3}) \left( \left( -(\varphi_{i+1} - \varphi_{i+2})^6 \right) - \left( -(\varphi_{i+1} + \varphi_{i+2})^6 \right) \right)\\
    =\frac{1}{3} (\varphi_i \cdot \varphi_{i+3}) \left( (\varphi_{i+1} + \varphi_{i+2})^6 - (\varphi_{i+1} - \varphi_{i+2})^6 \right)\\
\end{align*}
\begin{align*}
    \implies \frac{1}{3} (\varphi_i \cdot \varphi_{i+3})\left(\left( (\varphi_{i+1} + \varphi_{i+2})^3 - (\varphi_{i+1} - \varphi_{i+2})^3 \right)\cdot \left( (\varphi_{i+1} + \varphi_{i+2})^3 + (\varphi_{i+1} - \varphi_{i+2})^3 \right)\right) 
\end{align*}
using the difference and sum of cubes formulas

\begin{align*}
    \scriptstyle{ \frac{1}{3} (\varphi_i \cdot \varphi_{i+3})\left(\left(\varphi_{i+1} + \varphi_{i+2}-\varphi_{i+1} + \varphi_{i+2}\right)\cdot[(\varphi_{i+1} + \varphi_{i+2})^2+(\varphi_{i+1} + \varphi_{i+2})\cdot(\varphi_{i+1} - \varphi_{i+2})+(\varphi_{i+1} + \varphi_{i+2})^2]\right)}
\end{align*}
\begin{align*}
\\\scriptstyle{
\cdot\left(\left(\varphi_{i+1} + \varphi_{i+2}+\varphi_{i+1} - \varphi_{i+2}\right)\cdot[(\varphi_{i+1} + \varphi_{i+2})^2-(\varphi_{i+1} + \varphi_{i+2})\cdot(\varphi_{i+1} - \varphi_{i+2})+(\varphi_{i+1} - \varphi_{i+2})^2]\right)}
\end{align*}
\\
Simplifying

\begin{align*}
     \frac{2^2}{3}\left( \varphi_i \cdot \varphi_{i+3}\cdot\varphi_{i+2} \cdot\varphi_{i+1} \right)\cdot[(\varphi_{i+1} + \varphi_{i+2})^2+(\varphi_{i+1} + \varphi_{i+2})\cdot(\varphi_{i+1} - \varphi_{i+2})+(\varphi_i - \varphi_{i+2})^2]
\end{align*}
\begin{align*}
\cdot[(\varphi_{i+1} + \varphi_{i+2})^2-(\varphi_{i+1} + \varphi_{i+2})\cdot(\varphi_{i+1} - \varphi_{i+2})+(\varphi_{i+1} - \varphi_{i+2})^2]
\end{align*}\\

The integral reduces to factors, of which, by the deduction made in *(1), some of them must be multiples of 3. This proves that $P1$ is an integer. Now, since it has been proven that $P1$, $P2$, $P3 \in \mathbb{Z}$.

\begin{align*}
    \implies \int\limits_{\theta_2}^{\theta_1} f(x) \, dx  \quad \in \mathbb{Z}
\end{align*}

Similarly for $g$
\begin{align*}
=  \int\limits_{\phi_2}^{\phi_1} g(x) \, dx 
\end{align*}
    \begin{align*}
    =\int\limits_{\phi_2}^{\phi_1} 2 \cdot (\varphi_{i+1} \cdot \varphi_{i+2})x^2 + 2^2 \cdot (\varphi_{i+1} \cdot \varphi_{i+2}) \cdot (\varphi_{i+1}^2 + \varphi_{i+2}^2)x 
    + 2^3 \cdot (\varphi_{i+1} \cdot \varphi_{i+2})^3 \,dx
\end{align*}

Separating the integral 
\begin{align*}
    \int\limits_{\phi_2}^{\phi_1} 2 \cdot (\varphi_{i+1} \cdot \varphi_{i+2})x^2\,dx +\int\limits_{\phi_2}^{\phi_1} 2^2 \cdot (\varphi_{i+1} \cdot \varphi_{i+2}) \cdot (\varphi_{i+1}^2 + \varphi_{i+2}^2)x\,dx+\int\limits_{\phi_2}^{\phi_1}2^3 \cdot (\varphi_{i+1} \cdot \varphi_{i+2})^3 \,dx
\end{align*}

Defining for ease of notation $T1+T2+T3$ as each integral in order, just like in the previous proof, $T2$ and $T3$ turn out to be integers, since \\

\begin{align*}
   T2= \left[ 2 \cdot (\varphi_{i+1} \cdot \varphi_{i+2}) \cdot (\varphi_{i+1}^2 + \varphi_{i+2}^2)x^2\right]_{\phi_2}^{\phi_1}\\
   T3=\left[2^3 \cdot (\varphi_{i+1} \cdot \varphi_{i+2})^3 x\right]_{\phi_2}^{\phi_1}
\end{align*}

Integrating $T1$:

\begin{align*}
    \int\limits_{\phi_2}^{\phi_1} 2 \cdot (\varphi_{i+1} \cdot \varphi_{i+2})x^2\,dx=\left[\frac{2}{3} \cdot (\varphi_{i+1} \cdot \varphi_{i+2}) x^3\right]_{\phi_2}^{\phi_1}
\end{align*}

Applying the difference of cubes formula, we get:
\begin{align*} 
 =\frac{2}{3} \cdot \left( \varphi_{i+1} \cdot \varphi_{i+2} \right) \cdot\left( \left( -\left( \varphi_{i+1}^2 + \varphi_{i+2}^2 \right) + \varphi_i \cdot \varphi_{i+3} \right)^3 - \left( -\left( \varphi_{i+1}^2 + \varphi_{i+2}^2 \right) - \varphi_i \cdot \varphi_{i+3} \right)^3\right)\\
 =
 \frac{2}{3} \cdot \left( \varphi_{i+1} \cdot \varphi_{i+2} \right)  \cdot \left( -\left( \varphi_{i+1}^2 + \varphi_{i+2}^2 \right) + \varphi_i \cdot \varphi_{i+3}  + \left( \varphi_{i+1}^2 + \varphi_{i+2}^2 \right) + \varphi_i \cdot \varphi_{i+3} \right)\cdot
 \end{align*}
 \\
 \begin{align*}
     \left( -\left( \varphi_{i+1}^2 + \varphi_{i+2}^2 \right) + \varphi_{i} \cdot \varphi_{i+3} \right)^2+   \left( -\left( \varphi_{i+1}^2 + \varphi_{i+2}^2 \right) + \varphi_{i} \cdot \varphi_{i+3}\right)
 \end{align*}
 
 \begin{align*}
       \cdot \left( -\left( \varphi_{i+1}^2 + \varphi_{i+2}^2 \right) - \varphi_{i} \cdot \varphi_{i+3} \right) + \left( -\left( \varphi_{i+1}^2 + \varphi_{i+2}^2 \right) - \varphi_{i} \cdot \varphi_{i+3} \right)^2
 \end{align*}
\\

Simplifying the expression:

 \begin{align*} 
 =
 \frac{2^2}{3} \cdot \left( \varphi_{i+1} \cdot \varphi_{i+2}      \cdot\varphi_i \
 \cdot \varphi_{i+3} \right)\cdot
 \end{align*}
 \\
 \begin{align*}
    [ \left( -\left( \varphi_{i+1}^2 + \varphi_{i+2}^2 \right) + \varphi_{i} \cdot \varphi_{i+3} \right)^2+   \left( -\left( \varphi_{i+1}^2 + \varphi_{i+2}^2 \right) + \varphi_{i} \cdot \varphi_{i+3}\right)
 \end{align*}
 
 \begin{align*}
       \cdot \left( -\left( \varphi_{i+1}^2 + \varphi_{i+2}^2 \right) - \varphi_{i} \cdot \varphi_{i+3} \right) + \left( -\left( \varphi_{i+1}^2 + \varphi_{i+2}^2 \right) - \varphi_{i} \cdot \varphi_{i+3} \right)^2]
 \end{align*}

Similarly, for $g$, by deduction *(1), it is concluded that both $T1$, $T2$, and $T3$ have integer values, with $T1$ being an integer since some of its factors are multiples of 3, thus demonstrating that.\\
 
\begin{align*}
  \implies  \int\limits_{\theta_2}^{\theta_1} f(x) \, dx  \quad, \quad \int\limits_{\phi_2}^{\phi_1} g(x) \, dx\quad  \in  \mathbb{Z}\quad \blacksquare
\end{align*}

\section*{References}

\noindent
\begin{minipage}{\linewidth}
\small{
   $[1]$  Renault, M. (1996). The Fibonacci Sequence Under Various Moduli (Wake Forest University, Ed.) [Review of The Fibonacci Sequence Under Various Moduli]. Fibthesis.pdf; Marc Renault. \\
     \url{https://webspace.ship.edu/msrenault/fibonacci/fibthesis.pdf}
}
\end{minipage}

\vspace{0.5cm}

\noindent
\begin{minipage}{\linewidth}
\small{
   $[2]$ Aguilera, C. (2022). On Primitive Pythagorean Triples. ffhal-03247261v2. \\
    \url{https://hal.science/hal-03247261v2/file/Colorin_Primitive_Pythagoren_Triples.pdf}
}
\end{minipage}
\vspace{0.5cm}

\noindent
\begin{minipage}{\linewidth}
\small{
   $[3]$ Mathcentre. (2009). Quadratic Equations. mc-TY-quadeqns-1. mathcentre. \\
    \url{https://www.mathcentre.ac.uk/resources/uploaded/mc-ty-quadeqns-1.pdf}
}
\end{minipage}
\vspace{0.5cm}

\noindent
\begin{minipage}{\linewidth}
\small{
   $[4]$ Bicknell-Johnson, M. (n.d.). Pythagorean Triples Containing Fibonacci Numbers: Solutions for $F* \pm F* = K^2$. A.C. Wilcox High School, Santa Clara, CA 95051. \\
    \url{https://www.fq.math.ca/Scanned/17-1/bicknell.pdf}
}
\end{minipage}
\end{document}